\documentclass{article}
\usepackage{amsfonts}
\usepackage[polutonikogreek, english]{babel}
\usepackage[T1]{fontenc}
\usepackage[utf8]{inputenc}
\usepackage{amstext}
\usepackage{amssymb}
\usepackage{amsthm}
\usepackage{mathrsfs}
\usepackage{enumerate}
\usepackage[round]{natbib}
\usepackage{color}
\usepackage{graphicx}
\usepackage{graphics}
\usepackage{amsfonts}
\usepackage{amssymb}
\usepackage{amsthm}
\usepackage{mathrsfs}
\usepackage{amsfonts}
\usepackage{color}
\usepackage{amsmath}
\usepackage{latexsym}
\usepackage{graphics}
\usepackage{graphicx}
\usepackage{wrapfig}
\usepackage{cancel}
\usepackage{bussproofs}
\usepackage{amsfonts}
\usepackage{amssymb}
\usepackage{amsthm}
\usepackage{mathrsfs}
\usepackage{amsfonts}
\usepackage{enumerate}
\usepackage{color}
\usepackage{graphicx}
\usepackage{graphics}
\usepackage{amsfonts}
\usepackage{amssymb}
\usepackage{amsthm}
\usepackage{mathrsfs}
\usepackage{amsfonts}
\usepackage{enumerate}
\usepackage{color}
\usepackage{amsmath}
\usepackage{latexsym}
\usepackage{graphics}
\usepackage{graphicx}
\usepackage{wrapfig}
\usepackage{natbib}

\input{tcilatex}
\begin{document}

\title{Platonism, \textit{De Re}, and \\ (Philosophy of) Mathematical Practice}
\author{Marco Panza\thanks{Thanks to Jeevan Acharya, Annalisa Coliva,  Silvia De Toffoli, Michele Friend, Giovanna Giardina, Peter Jipsen, Alex Kurz, Abel Lassalle Casanave, Fabio Minazzi, Andrew Moshier, 
Michael Otte, Andrea Sereni, Stewart Shapiro, Daniele Struppa, Ed Zalta, two anonymous referees, and countless  other friends and colleagues.} \\ CNRS, IHPST (CNRS and University of Paris 1, Panthéon Sorbonne) \\ Chapman University}\medskip
\date{\today}
\maketitle

\section{Preamble}
Mentioning platonism and \textit{de re} (attitude) in the title of a chapter of an
handbook about philosophy of mathematical practice might surprise some adepts of the latter. While platonism 
might be seen as a metaphysical option, the distinction between \textit{de re} and \textit{de dicto} attitudes 
might appear as proper to the analytical philosophy of language. These to facts might be considered enough 
for making the philosophy of mathematical practice principally alien to these notions. 
In echoing this view, while reporting on a previous version of my essay an anonymous referee 
has qualified it as ``an ultra-conservative piece in the philosophy of mathematical practice''.
The reason for this is that it would have told ``nothing about mathematical practice, divorcing itself starkly 
from the aspirations of some philosophy-of-mathematical-practice pioneers such as the PhiMSAMP \citep{PhiMSAMP} 
members, nor does it agree with Mancosu’s (\citeyear{Mancosu2008})  conservative 
understanding of PMP as exploration of new philosophical questions about mathematics''.

This report made me decide to write the present preamble, in whose absence  some readers might 
have been tempted to agree with this judgement, by misconstruing, then, the very purpose of my essay.

It is certainly true that I hereby explore a traditional (or ``ultra''-traditional) question for the philosophy of mathematics: 
so traditional to date back
to ancient Greek philosophy. However, being new is a required value neither for good philosophical
questions nor for fruitful philosophical approaches. A philosophical investigation is 
often made valuable and innovative because of the way traditional questions are re-stated and tackled
and of the tentative answers that stem from time-honored approaches. It is certainly not to me to say whether
the investigation I present here is a valuable and innovative one. What I can say is that it aspires
to be so because of this reason: since it aims at using a well-established distinction, as that between
\textit{de re} and \textit{de dicto} attitudes, to answer a traditional question in a purportedly original
way.

It is also certainly true that I hereby explore no specific aspect or episode of mathematical practice. Promoting
and arising similar explorations is not, however, the only purpose of philosophy of
mathematical practice (as I understand it, at least). Another equally important purpose is accounting for some
essential general feature of this practice, by so insisting that the term `mathematics' in `philosophy of
mathematics' should stand for it: for the human activity of doing mathematics. This activity has
social and material aspects, but it is not because of them that it is
mathematical. What makes it so are the modalities of its being an intellectual activity. This makes it necessary for the
philosophy of mathematical practice necessary to engage with the purpose of accounting for the specific feature of that
intellectual activity that we perform by doing mathematics, namely with the question of wondering what makes
an activity mathematical in nature. This is just the question I explore in my essay (even if only partially, of course). And it
is just because of such that I take this essay to be, in the full sense of the term, a piece of philosophy of mathematical practice—whether 
(ultra-)conservative or not, I leave to the readers to decide.

More than having a history, mathematical activity is (or generates) a history. Since I in no way think that
the philosophy of mathematics should have any normative duty concerning mathematics itself, I take this history
to be, properly speaking, what the term `mathematics' should stand for in `philosophy of mathematics'.
 Establishing what and
how mathematics should be is a question for professional mathematicians, not for philosophers. 
Even less, it is a duty of philosophers to establish whether something is there or not.
In particular, philosophers of mathematics are entitled to decide 
neither what is actually there (or has happened and happens) 
that one would consider as mathematical, nor 
what is genuinely mathematical and what is not so, among
all that is there (or has happened and happens) that  mathematicians and professional historian of mathematics call by
this name.

In my view, philosophy essentially aims at accounting for
independent realities (unless it is some sort of meta-philosophy, which it might of course legitimately be, in some cases). 
This is even obvious, it seems to me, when it is the philosophy of something, as the philosophy of mathematics
clearly is. Since, in order for philosophy to be the philosophy of something, this something should be available to its
inquiry.  This inquiry can (or even cannot avoid to) shape this something, just as science do for the reality it studies.
But this shaping is just an essential part of accounting for it, of explaining how it is, of spelling it out, not of establishing 
how this should be, or what is there to form this reality, and what is not there. 

Arguing from these theses is, in my view, just the same as promoting a practical turn in philosophy, if there is one that one
might promote. 
Hence, if philosophy of mathematical practice is to be part of this turn, it has to tell us
of which something that is there and is called `mathematics' it is intended to be philosophy of. When I say that
the term `mathematics'  in `philosophy of mathematics' should stand for the historical fact of mathematical activity,
I simply mean that this something should be nothing but this historical fact: what this activity has been and is today. 
This fact makes the intersection of
history and philosophy of mathematical practice not only quite strict, but essentially constitutive: the philosophy of
mathematical practice is, to my mind, the philosophy of an historical fact.

Still, there is no need, I hope, to refer
to Croce (the historicist philosopher \textit{ex excellentia}) for remarking that history (\textit{res gest{\ae}}) is not the same 
as historiography (\textit{historia rerum gestarum}).  
This constitutive connection does not make philosophy of mathematical practice the same thing as the history of mathematics,
in the usual sense where `history' wrongly stands for `historiography' (and, then, philosophers of mathematics the same 
thing as historians of it). History (in this sense, i.e. historiography) reconstructs
the relevant facts and accounts for their specific and detailed aspects. Philosophy accounts for general features,
the way of being, the peculiarity that makes them mathematical facts. The former can be done without the latter
(though possibly not very lucidly). The latter cannot be done without the former, since the former provides to the
latter with the reality which it is the philosophy of (since history, properly called, only appears to us through the lens of historiography:
we have no direct access to it). This is another way to say, in the the language of Lakatos's dictum, 
that historiography of mathematics without philosophy of mathematics is certainly blind, but can still subsist, while
the philosophy of mathematics without the history of mathematics (and therefore its historiography that provides it to us) is
so empty that it cannot subsist unless as a caricature of itself.

Another important clarification is suggested by the same referee of 
a previous version of my chapter. What I'm arguing for is not at all that the philosophy of mathematics cannot but
consist of a collection (or even an organic system) or case studies. My point is precisely the contrary. 
A collection (or system) of case studies cannot but be what it is: namely a more or less ordered configuration 
of particular inquiries devoted to single neighborhoods in the space of history. Philosophy is
much more and much less than this: it is an account of the space itself. While the
former accounts for some elements in a quite large set, the latter defines a structure based on this set
and studies it as such. This is why the reader will find no case study in what follows (while some of them can be easily
found in other, more properly historical essays I have written), but rather a reflexion of what makes
the historical fact of mathematical activity mathematical in nature.

I conclude my preamble with a proviso. Though I'll be naturally brought to mention
a few philosophers of mathematics and some of their works,
my purpose here is not that of  discussing the opinion of this or that 
colleague or master. Rather it is that of promoting a general philosophical option, and
rejecting some others,
by presenting the former in  my personal (and I hope, original) guise. Any reader that would like to have 
more bibliographical information about the present discussion in the philosophy of mathematics that
I shall deal with should look for it elsewhere. I contributed to offering a survey of some of  these discussions in 
\citep{PanzaSereniPlato}, and (apart fort the most recent developments, elapsed after its publication) I
can certainly refer the reader to this work.

\section{A Misconception\label{Sect.Mis}} 

\noindent Philosophy goes with many temptations. And the philosophy of
mathematics is not without them, too. The most dangerous one is taking the
latter as a mere extension of the former, due to the possibility of using
mathematics as a testbench  for philosophical theses, as a storehouse of
examples and problems on which to measure the mutual merits of opposite
philosophical views. There are two kinds of philosophers: those that have
problems and look for solutions for them, and those that have solutions and
look for problems to solve them with. The philosophers of the latter kind
are victim to a misconception, I contend.
The dangerous temptation I'm alluding to is a form of this misconception. One of the ways it
takes place in the philosophy of mathematics is using mathematics and/or
mathematical examples to measure the mutual merits of realism and
nominalism, understood as general metaphysical theses, the former asserting
and the latter denying the existence of abstract objects.

In the course of this exercise, realism, or, more precisely, the
particular version of it that is often called `platonism'---which I write, as it often done
(following \citealp{Bernays1964}, 1964, p.~275), with
a lowercase `p', for disguising it from other theses more faithful,
also in their intention, to Plato's philosophy---is presented as
the thesis that ``there exist mathematical objects'' \citep[p.~1]{Balaguer1998}.
More precisely, it is the thesis that, independently of our course of action and thought,
there exist abstract objects in their own right, which are intrinsically
mathematical in nature, and that mathematics deals with them, and aims at 
describing them, by so revealing how they are and how they are related 
to each other.
Mathematics would, then, be successful when what it asserts of these
objects is true, that is, when it truly describes them and their world as
they are, in fact. It follows that a mathematical theory is good if and only if its
theorems are true of the relevant objects. Among these objects there are
sets, numbers, and geometrical figures, or, possibly, only sets, starting
with the empty one, according to the most reductionist ones among these
realist self-styled philosophers of mathematics. 

These, or other, less
radical, partisans of ontological reductionism are certainly ready to admit
that good mathematical theories can include theorems that are not literally
true, or are so only if the relevant terms are properly understood according
to appropriate reductionist clauses. Still, once the reduction is finally
complete, or these clauses are fully operative, there is no more need for
rephrasing: the theorems of these theories plainly means what they appear to
mean \textit{prima facie}, and what makes them true of the relevant objects
is that these objects just are and stand to each other like these theorems literally say that they are.

Nominalism is presented as the opposite thesis that denies the existence of
abstract objects, and, \textit{a fortiori}, of mathematical ones. It denies, for instance,
 ``that numbers, functions, sets, or any
similar entities exist''  and, by consequence, ``that it is
legitimate to use terms that purport to refer'' to them \citep[p.~1]{Field1980}.
To avoid the
most classical of paradoxes, this cannot but be understood as the thesis
that the concept of an abstract object is
empty, which immediately entails that this is also so for the concept of a
mathematical object (provided it be admitted that mathematical objects could not 
be but abstract ones). These two concepts 
cannot be taken as ill-formed of course, since,
if they were so, it would be hard to sensibly claim that they are empty. The same is even more clearly
the case for more specific mathematical concepts, as that of a number,
a set, or a function. These concepts are so defined to make 
some particular mathematical items purportedly fall under them.
It follows that either a nominalist is prepared to also argue that these
items are not objects (but, then what?), or cannot avoid to conclude that
the usual mathematical
language is both well-formed and empty, making mathematical
statements either literally false, or only vacuously true. 

Supposing that the doubtful distinction between mathematical items and objects is 
drop out, it follows that appropriately understanding a mathematical statement depends
on two possible moves, locally alternative to each another, but
possibly operating together within a more general program aiming at a
global understanding of mathematical discourse. On the one side, one
might require, in a Quinlan vein (inaugurated in \citealp{Quine1939},
but the pursued in many subsequent works), that mathematical statements be paraphrased
as statements on non-abstract, and then
non-mathematical, objects. Once
the paraphrase is complete, these statements would become, finally, literally true
or false. This would allow one to say that a mathematical theory is good if its
theorems admit a paraphrase making them literally true, but no more
mathematical, statements. One might also require, on the other side, as suggested in
\cite{Field1980},
that
mathematical theories be added to empirical ones, which only deals with concrete
objects to generate conservative extensions of them. Though remaining
literally untrue and not necessarily admitting any paraphrase
making them true, the theorems of a good mathematical theorems
would then have a very welcome virtue: that of making it easier (or even
practically possible) to establish empirical truths.

Both the partisans of platonism and nominalism so conceived  may easily find, within the mathematical
storehouse, many examples appearing, at first sight, favorable to one or the
other party.

A nominalist might easily observe, for example, that a sentence like $%
\ulcorner$Mars has two moons$\urcorner $ can be faithfully rephrased as $%
\ulcorner$There are $x$ and $y$ such that they are distinct to each other and both moons of
Mars, and for any $z$, either it is not a moon of Mars or it such that
either $z=x$ or $z=y\urcorner $, where no reference is made to any number.
The use of the term `two' might, then, be taken as a \textit{fa\c{c}on de
parler}, simply used for short. 

A platonist might concede this, after all,
but observe that the former sentence does not express, in fact, a (purely)
mathematical statement, and that the simplicity of such an example vanishes,
when we speak of numbers as such, rather than of the numbers of moons of
Mars, or of pears or apples. When we say, for instance, that $2$ plus $3$ is $5$, we
do not speak at all of two particular things and three other particular things that, when put together,
result in five things. At most we speak, in agreement with def.~VII.2 of Euclid's \textit{Elements} 
(\nocite{EuclidH}Euclid, TBEH),
of two units and three units, that, when put together, result in
five units. But units are certainly not, as such, in all their
generality, concrete objects. This makes it impossible---the platonist might conclude---to 
follow the same strategy used above to restate the sentence to avoid any reference 
to abstract objects. Indeed, this would make it loose  its generality, which is a salient feature of 
mathematical statements.

The nominalist can try to replay by
suggesting a new paraphrase that would be too long to write here. It will be enough
to say that it might try to render the generality of the original sentence by
replacing the reference to units either by appropriate universal quantifiers ranging
on the totality of concrete objets, or by a talk of concrete objects whatsoever. At this point, the platonist
might either present another simple example, as that of the sentence  $\ulcorner$there
are infinite prime numbers$\urcorner$, by challenging the nominalist to
suggest a faithful nominalist paraphrases, or, go for a much more general
strategy and
remember Plato's distinction between vulgar and knowledgeable (or
philosophical) arithmetic, respectively concerned with ``unequal units 
[{\greektext{μόνᾳ ἄνισος}}]'', such as ``two armies'' or ``two oxen'', 
and with ``countless units [{\greektext{μόνᾳ μυρίαι}}]'', none of which
differs from each other (\textit{Philebus}, 56d-e;
see also \textit{Theaetetus}, 195e-196a). The nominalist might retort, in turn, under
Quine's authority \citep[p.~400]{Quine1986}, that the latter arithmetic is merely \textquotedblleft
recreational\textquotedblright and should in fact be ignored by
naturalized (and, then, good) philosophy, to the effect that an
example as the last one just mentioned, should be taken as perfectly immaterial for
such a (good) philosophy. 

This last move might also open the
road for the second nominalist strategy, in which the recourse to paraphrase is
replaced by the talk of conservative extensions of scientific, and, then,
nominalistic, theories. Faced with this strategy, the platonist would have no reason to
surrender, since there would be still room to argue: \textit{i})~that scientific
mathematized theories are far from nominalist, in fact, since mathematics is
constitutive of its objects; and \textit{ii})~that pure mathematic is far from
recreational and should constitute the appropriate text-bench for
deciding the realism vs. nominalism dispute.

Who is right? It might appear that there are good reasons for both sides.

For instance, there is no need to invoke the number two, as an abstract
object, to count the moons of Mars if counting objects (or transitive
counting: \citealp{Benacerraf1965}, pp.~49-50) is understood as such, that is as an intellectual operation on
these objects (rather than as the establishment of a bijection between them
and the numbers), which is all that is needed for all scientific purposes. 
There is no more any need to invoke numbers to decide whether there are more
planets or moons in the solar system. Using them is surely quite
convenient; it makes things simpler. However, it is not indispensable. If one had
to work on much bigger collections, one would probably need to pass through
numbers, but this would be only a practical need. In principle one could
avoid them. This might appear to be a strong enough reason to think
that number words enter our statements concerned
with these empirical practices only as useful, but dispensable, 
\textit{fa\c{c}ons de parler}, as the nominalist suggest. 

But is it the use in
such empirical practice that decides of the nature of mathematics? For
making the realist win over the nominalist on the testbench  of mathematics,
it does not seem at all necessary that no occurrence of mathematical words might be
explained as the nominalist suggests. It seems enough that some
occurrence cannot be so explained. But then, the case of the infinity of
prime numbers might become decisive: either the nominalist can
explain it nominalistically, or the platonist has a point, after all, unless, in a last and
decisive gasp, the nominalist succeeded in excluding this case from the relevant
testbench.

\section{An Ill-Posed Question}

But, then, again: who is right? Neither one, I contend. Since the question
the two parties offer opposite answers to is simply ill-posed. 

Even supposing that
the realist vs nominalist dispute has any sense, in general, it does 
not extend to mathematics in these terms, at least if such an
extension has to take the form of a philosophy of mathematics. Since, if
philosophy of mathematics is to be about mathematics, it has to start from
mathematics itself, and reflect on it, rather than begin elsewhere and come to
mathematics by accident, by reproducing on it the same general
questions and answers already looked over outside it, 
with the mere replacement of schematic place holders with mathematical
ingredients. In other terms: it is not enough to replace `abstract' with
`mathematical' in the general question `do abstract objects exist?', to get
a sensible question for the philosophy of mathematics.

Let us be a little bit more precise. Let us assume that
the word `abstract' is the only schematic word in this general
question. This means that the replacement of
this word with any other word of the same type within this question is not to
go together with a
corresponding replacement of any other words such as `objects' or
`exist', or, at least, with a corresponding particular understanding of
them. This being granted, suppose, for the sake of
the argument (and only for that), that, when it is so understood, such
a general question has a clear sense. The point is then the following: also
under these provisos, it would still remain that the 
mathematical instance of the question---`do mathematical
objects exist?'---is not a genuine question for philosophy of mathematics, or,
at least, it is not so if it is also admitted that philosophy of mathematics is actually 
to be about mathematics. 

Of course, philosophers are perfectly allowed to have no particular interest for
mathematics, as well as for any other field of knowledge.
But we, philosophers of mathematics, are also allowed, and to my mind, even required,
to regard the philosophy we practice as a genuine enterprise only if it reasons on
mathematics as an external material that it takes as given as such, before beginning its
work.

If this is admitted, there is room to argue that the very question `do
mathematical objects exist?' is nonsensical. Or, at least, it is so, 
if mathematics is
supposed to be given, and it is granted that the words `objects' and `exist' 
designate primitive notions, asking for no appropriate elucidation depending on 
mathematics itself. Consequently, the connected question `is
mathematics about its own objects?' also becomes nonsensical, under similar
conditions. 

This makes any putative answer to these questions, advanced in
agreement to such conditions, not only doubtful, but inevitably misleading: a
claim that is not so far, after all,  at least for its form, from Kant's claim about
dogmatic metaphysics.

A symptom of this nonsensicalness is the difficulty of answering (and 
in agreeing on the good answer) to a connected question: `what could be
conceived of as an evidence for the existence, or non existence, of
mathematical objects?'. If mathematics is taken as given, it should include
such evidence. Of course, this might be very well hidden,
quite far from the surface that immediately appears to an
inexperienced observer.  This would make discovering and exhibiting such evidence quite
hard. But it should be, in any case, possible to
clearly say how it might be like. Since the question `do mathematical
objects exist?' presents itself, as a factual question, deciding which is
the right answer to it
should be a matter of fact. Moreover, assigning a clear meaning to this question, by making it sensible, should
require clearifying how this matter of fact could become manifest to our
understanding. We should, in other terms, be able to say what would make us
unquestionably aware that mathematical objects do or do not exist.

I know
only two attitudes face to this last crucial question, which might appear plausible, at last at first glance. 

The first is only available to nominalists: the burden of providing evidence for
deciding about the existence of mathematical objects is only on the platonist side, since platonists 
are the only ones
that make a positive claim.  Here is how the point is made by 
C. Cheyne (\citeyear{Cheyne2001}, p.~108) for instance:
\begin{quote}
[\ldots] it is a contingent matter whether numbers exist with mathematical
necessity in the actual world. The claim that this could be decided by
purely a priori methods becomes dubious. Discovering whether or not this
world is favored by the existence of entities that might not exist would
seem to be a matter of observation. But how could we observe such entities
if they lack causal powers? The burden of proof is [\ldots] on platonists [\ldots].
\end{quote}

The nominalist's burden is then nothing but that of accounting 
for the available evidence in a way that does not suppose the existence of an abstract,
and, then, mathematical object. This is not properly an answer to the
question; it is instead a move for turning the absence of any answer in favor
of the nominalist side. But the move is weak enough. For it turns this
absence in favor of an agnostic attitude, at most, by so making clear that
the nominalist thesis is, in fact, a \textit{petitio principi}. The nominalist
might retort by appealing to the Ockham razor, and claim that, in the absence of
any positive evidence, all that can decide whether mathematical objects exists or not 
is whether it is necessary to postulate their existence in order to account for the
available evidence. However, this inverses the burden of proof, since it now requires
the nominalist shows that this is not necessary. In J. Burgess's words \citep[pp.~95-96]{Burgess1983}:
\begin{quote}
Actually, the burden of proof is on such enemies of numbers [\ldots], 
to show either: (\textit{a}) that science, properly interpreted, already does
dispense with mathematical objects, or (\textit{b}) that there are scientific reasons
why current scientific theories \text{should be} replaced by alternatives dispensing
with mathematical objects.
\end{quote}

The problem becomes, then, that the nominalist account, both in terms or paraphrases and
in terms of conservative extensions, is partial, at most, if not largely
unsatisfactory. Possibly, it accounts for some scant evidence. But it does it
at the price of declaring this evidence untrustworthy, that is, of tacking it
as a fake appearance of a hidden reality that our language is often
incapable to describe faithfully, or can only describe in a quite entangled and
odd way.

The second attitude is that of some platonists: the required (positive)
evidence is provided by the truth of some appropriate mathematical statements
involving terms that purportedly denote mathematical objects. in
B. Hale and C. Wright's words \cite[p.~8]{HaleWright2001}: 
\begin{quote}
[\ldots] no more is to be required  in order for there be
a strong prima-facie case that a class of apparent singular term have reference,
than they occur in true statements free of all epistemic, modal, quotational, and
other forms of vocabulary standardly recognized to compromise straightforward
referential function.
\end{quote}
Since a statement about
certain objects is true only if these objects exist. Fair enough.
But this leaves other questions open: what evidence do we have for the
truths of the relevant mathematical statements? What evidence 
do we have that these statements  are
actually about the objects they seem \textit{prima facie} to be about? Even more
crucially: are we
really sure that mathematical truths require existence? 

The second of these three questions
brings grist to the mill of nominalists. Considering the third one would
bring us too far from the purpose of the present paper. So, I limit myself
to the first. A possible answer to it is that the relevant truths are
logical truths, or, at least, logical consequences of so innocent
assumptions that we could not reject them as untrue. This was essentially
Frege's strategy, and, in a different version, is still the neo-logicist
one (as defended in \citealp{HaleWright2001}, for instance). 
But it has problems, in turn, and not only because of the difficulty of
deriving suitable mathematical statements from consistent systems of innocent
enough assumptions. Also, among other things, because of the doubt one can
raise about the fact that the innocence of these assumptions 
entails their truth, or, better,  an appropriate
kind of truth for supporting existence. One might, of
course, redefine existence in terms of (the relevant kind of) truth. But this
would violate the condition mentioned above, about the
presence of a single schematic word in the question about the existence of
abstract objects, making the platonist thesis
essentially different from a simply application to mathematics
of the realist one. Under this guise, the thesis would rather involve
a claim about a peculiar notion of existence, proper to mathematics. It
should, then, be supported, not by a general metaphysical attitude, but
by a specific inquiry about the nature of mathematics, which is exactly what I'm encouraging.

\section{Plato's Problem}
I can stop the \textit{pars destruens} of my paper here. It was aimed at critically
describing an attitude that no sort of philosophy of mathematical practice might share,
to my mind. Still, arguing that the philosophy of mathematical practice 
might not share this attitude 
is not the same as arguing that the question it deals with, though nonsensical, as such,
is not rooted in a relevant and perfectly sensible an genuine problem about mathematics (which
philosophy of mathematical practice is certainly not only entitled, but also required  to tackle).
What I wanted to suggest is that the trouble with this attitude is much less due
to the opposite answers that it suggests for this question, than to the question itself, better to the way 
this question renders the crucial problem that stays at the roots of it. 
The point at issue is that this way makes this question 
essentially external to mathematics, by so making hard to
decide which evidence, coming from mathematics itself might count as a
justification for an answer to it, whatever this answer might be.

In \citep{PanzaSereniPlato}, I designated the problem I'm referring to as Plato's
one. A way to address it is by quoting a famous passage by the 
\textit{Republic} (527a-b), in an appropriate translation:

\begin{quote}
This <\ldots> will not be disputed, at least, by those that <are> even a bit
acquainted with geometry, that this science [{\greektext{αὕτη ἡ ἐπιστήμη}}] 
is in full opposition [{%
\greektext{πᾶν τοὐναντίον}}] with the language spoken within it by its
adepts <\ldots>. Since <the way they> speak <is> as much ludicrous as
necessary [{%
\greektext{μάλα γελοίως τε καὶ ἀναγκαίως}}]: they shape all their discourses as 
if they were operating in practice and
in the purpose of practice [{%
\greektext{ὡς γὰρ πράττοντές τε καὶ πράξεως ἕνεκα πάντας τοὺς λόγους ποιούμενοι}}]; 
[they] speak of squaring, applying and adding, and they all say this way, 
whereas all <their> learning is somehow suited for the purpose of knowledge [{%
\greektext{τὸ δ᾽ ἔστι που πᾶν τὸ μάθημα γνώσεως ἕνεκα ἐπιτηδευόμενον}}]. <\ldots > [This is] knowledge
about that which always is, and not that which now comes into being
and now perishes.
\end{quote}

Contrary to the way this passage has been often understood 
(and consequently translated), this is in no
way a criticism of the language and practice of geometers (the point
has been firstly made, \textit{mutatis mutandis} by M. F. Burnyeat, \citeyear{Burnyeat1987}). 
Plato notices
that this language seems to be shaped ``{\greektext{πράξεως ἕνεκα}}
[in the purpose of practice]'', whereas geometrical learning is ``{\greektext{γνώσεως ἕνεκα}}
[for the purpose of (theoretical) knowledge]'' of eternal (and abstract)
matters. He takes this language to be``ludicrous [{\greektext{γελοίως}}]''. But also 
``necessary [{\greektext{ἀναγκαίως}}]''. He
seems then to acknowledge that there is no other option than using such a practical,
ludicrous language for theorizing about such an abstract matter. 

This is exactly the
opposite of the nominalist claim: mathematics does not use an abstract
language for speaking of concrete objects, by calling for appropriate
nominalist paraphrases; it rather speaks of \textit{abstracta} as if they
were material artifacts, by so calling for an appropriate understanding. The
ludicrous (or metaphorical) aspects of its language does not consists in
its surreptitiously dealing with concrete objects by \textit{prima facie}
speaking of abstract ones, but in its (having no other option than) using
a language originally fashioned to speak of the former objects for reasoning 
about its abstract matter. 

But this is no more the same as plainly asserting that
there are abstract objects that mathematics is about and call with their proper
names. If it were so, there would be no opposition, no necessary
ludicrousness. At most, the problem would be that of inventing a new language, just fashioned
for speaking of these abstract objects, just as the vulgar one had been
invented to speak of concrete ones. 

According to Plato, mathematics is certainly not about ``ideas [{\greektext{εἴδη}}]''.
It cannot depend on the exercise of {\greektext{νοῦς}}, but rather depends on that of {{\greektext{διάνοια}}, or,
possibly {\greektext{φαντασία}}, as Proclus will say later.
The soul can contemplate, at most, the idea of a triangle, or that of a number, possibly those
of even and odd, but mathematics requires more, or less than that: it requires triangles, and 
numbers; not single ideas, but pluralities of items. What do these pluralities consist of? How can we
conceive of these items? More in general: how can we account for the intellectual exercise that
mathematics pertains to? 

Such is Plato's problem. And this should be the ontological problem of the philosophy of mathematics. Not
whether there are {\greektext{μαθηματικά}} or not. But how can we account for our reasoning on them?
Since it is of them that mathematics is reasoning, unless one wanted to argue, \textit{contra}  Plato, and because of a
biased and pre-established misconception, that its form is fictitious and is to be amended by true philosophy (which would
just be the opposite of philosophy of mathematics, being rather philosophy against mathematics).

Take this statement $\ulcorner$any vector space has a basis$\urcorner$. What is it speaking of? The most obvious answer is
that it speaks of vector spaces and its bases. Does it not? Why not? Possibly since vector spaces must be abstract
objects and abstract objects do not exist? I suggest a nominalist supporting this answer to explain it to any 
working mathematician, and
consider the reaction. I suspect that it would be better not to describe it. But, does this mean that the platonist described above has finally
won? Does this prove that vector spaces exist and have bases? And that there is nothing more to say on that matter, except, possibly,
for adding that this makes the statement true, and its proof a justification for its truth? 

Suppose to work within \textsf{ZF}, i.e. whiteout the axiom
of choice. The proof that any vector space has a basis does no more hold. 
But what does the axiom of choice say? Let us consider one of its most usual formulation (all equivalent to
each other in appropriate frameworks): $\ulcorner$For any set of pairwise disjoint non-empty sets, 
there exists at least one set that contains one and only one element for each set in the former set$\urcorner$. What is it speaking of? Again,
the most obvious answer is
that it speaks of (non-empty) sets, or, more specifically, of sets of (non-empty) sets. Our platonist would say that these sets exist 
(let us avoid asking him/her whether this
makes the set exist formed by all sets). But, then, is this statement true of them? If it were not, I suppose the same platonist should 
maintain that it is false.
Then its negation should be true. But, based on this negation, we might now prove that there are vector spaces with no bases. 
Hence our previous statement
would be false. The former statement ($\ulcorner$any vector space has a basis$\urcorner$) is, then, true only if the latter (the axiom of choice) is so. 
But how might we prove the latter? We might certainly prove it, in \textsf{ZF}, by
grounding on another, equivalent or even stronger statement, for example on the former one (which is provably
equivalent to the former over \textsf{ZF}). But
this would only displace the question. What is it if it is not a proof that justifies the axiom of choice? And if we have no 
justification, how can we be sure that it is true. And if we cannot be sure that it is so, how can we be sure of it for the former statement, about vector spaces?
So, possibly the two statements are not true: possibly vector spaces exist, but some of them have no bases. All that we can know, then, 
according to our platonist, is that vector spaces exist and either all of them have a basis, or some do not. 
Or, possibly, this is not so, since, in fact, we have, after all, a clear justification of the axiom of choice: we see that sets are as it says; we see it with
our insightful eyes of mathematicians, even, if we cannot make other eyes see it, even eyes of other (non-classical or constructivist) mathematicians.}

The problem with both conclusions is that they do not follow from any evidence coming from mathematics itself. They do not really provide a solution
to Plato's problem (relatively to the particular fragment of mathematics that it is concerned with), but merely dissolve it. 
If these conclusions hold, the
problem would simply not arise; it would be solved before being posed. 
Since, what causes the problem is just that mathematics speaks of vector spaces just 
like geography speaks of lakes, mountains or towns, whereas vector spaces are not like lakes, mountains or towns. 
So, how can mathematics speak of them
this way? The previous platonist conclusions simply deny the root of the problem: they merely presume that  vector spaces are just like lakes, mountains or 
towns, at least at the insightful eyes of (classical) mathematicians. 

\section{On the Platonist Side \label{sect.PS}}

Let us restate the problem, then, in a more general setting. Mathematics results from a quite 
peculiar intellectual activity, or it is such an activity when it is considered as an enterprise rather than as a corpus. 
This activity appears to be twofold. It seems to consist both in the creation from noting of eternal and immutable  items, and in the
attribution of (some) properties and relations to them. 
Though these items appear to be eternal and immutable, their 
creation is the work of
human beings (not of Gods or God), and it takes place in time, or even in history 
(not before any time). These human beings have, moreover, no other way for justifying 
their attributions of properties and relations to these items
than proving statements that assert that these items are and stand
to each another accordingly to these attributions. 
The proofs of these statements start from other statements of the same sort,
some of which have necessarily to be admitted both without grounding on a
similar justification,  and without having acquaintance with these same items. 
How can we explain it?

This is not at all a problem about existence, about what exists or does not exist. What exists is clear: it is mathematics, and this existence is made manifest by 
history and actual experience (of teachers, students and professionals), just as it happens for the existence of other intellectual 
enterprises and social phenomena,
as well as for their products, and, \textit{mutatis mutandis} (but not \textit{mutandis} so much) for empirical middle-size objects. The problem is rather about
how this existence should be better accounted for.
As such, this is a genuine, even crucial problem for a philosophy of mathematical practice, since this
existence is just that of  mathematical practice itself and of its outcomes. 

When one looks at the problem this way, the theses that the platonist and nominalist I have described above try 
to test (or, better, confirm) on the testbench  of mathematics 
are invested in another light, and they appear, then, dramatically naive and implausible. The platonist thesis becomes the thesis that this evident reality is 
to be explained by postulating the exisrence of another reality, inaccessible and unverifiable, but somehow (either causally or not) operative 
enough so as to make it possible to
take the former existence as a more or less deformed imagine of the latter.  
The nominalist thesis becomes the thesis that this evident reality is not as it looks, since what appears is a
theater performance disguising a much more trivial reality, not dissimilar in nature from our everyday life. 
When a serious philosophy of mathematical practice tackles the
problem, it cannot but reject both these preposterous answers.

To my mind, the rejection should not be symmetric, however, since these theses are not so: whereas the former aims at explaining the available evidence through a
farfetched assumption, the latter aims at denying this evidence in the name of a preconception. If I had to take a side in the dispute, though judging
it misleading, I would have, then, no doubt in taking the platonist side. Since this is  the side which
aims at accounting for mathematics as it appears to be, rather than as
it is required to be in agreement to an independent preconception. The purpose of the following considerations is to suggest
a plausible way to stay on this side for
a philosophy of mathematical practice: a way of being a platonist without recourse to the  
farfetched assumption of the existence of mathematical objects.

One might argue that staying on the platonist side because of the reason I just advanced is still not the same as being a platonist. If it were so, the  
preceding pronouncement would be incongruous, since what comes after the colon would go far beyond what comes before it. Anyone that
agreed with this would also probably find the position I will defend in what follows far from platonist in nature. 
Though I do not care so much about -ist and -ism-words,
I regard this point as relevant, and I want to clarify why. 

To my mind, this is not a purely terminological point, indeed. It rather depends on a judgement on what is essential in
the platonist option, in all the different forms it  has taken along the history of philosophy of mathematics. More than 
that, it also depends, in a sense, on what is
considered to be the main task of philosophy. For many, this task is that of establishing what there is, and platonism is an important part of philosophy
since it clearly claims that there are abstract objects. Not for me! I do not take this to be at all a task for philosophy. At most it is one for positive sciences
(whether they be physical, natural,  biological, economical or historical). I neither take the task of philosophy 
to say how what is there should be (if it could be in different ways). This is certainly
a task for politics, but philosophy is not politics at all (though it could be useful for it). Philosophy should rather aim at providing the most 
appropriate conceptual categories to be used to account for what there is, or, better, what appears to there be.

This makes me think that platonism is a worthwhile option in philosophy of mathematics since it
suggests that mathematics deals with the very contents that the language of mathematics speaks of, \textit{prima facie}. In other terms, I take
platonism, in the philosophy of mathematics, to be the quite general thesis that mathematical talk is to be literally understood as a talk of individual contents that
provide the subjects of different form of predication---in a word, as a talk of objects---, and that this talk, so understood, consists of statements that hold
as such, and as such have to be admitted and accounted for. My purpose, in what follows, is to articulate and better detail this general thesis in an
appropriate way, along the line of a philosophy of mathematical practice.

\section{Philosophical Purposes}

Before doing that, an excursus is in order, however. It aims at shortly answering this further question: 
what should the purpose be, in general, of a philosophy of mathematical
practice? Or even, more simply: what should such a philosophy be like? I have already tackled the problem in my preamble, 
but it seems in order to come back on it now,
in the light of the previous considerations.

Let us put the question this way: should philosophy of mathematical practice be conceived as a reaction against a metaphysical, analytical, or foundational trend
in the philosophy of mathematics? Or rather as an attempt to tackle new, previously unnoticed problems that, though strictly connected with the everyday work of
mathematicians, do not strictly have a mathematical nature, being rather concerned with methodological choices, value judgments or alike? Or, again, as an effort to 
rethink classical problems of philosophy of mathematics, restate them, and tackle them from a different perspective? 
 
I do not like the first, purely antagonistic, option, and consider that the second one is simply non enough. 
Since I think that classical philosophical problems are so for they
are good and difficult problems, though they are often conceived and restated in an ill-posed way. Plato's problem is an example of that:
what makes it be, since Plato himself, a traditional problem for philosophy of mathematics, is that it is a crucial problem that no serious philosophy of mathematics can 
dismiss, and even less if the dismissal is simply motivated by its being an old problem. This makes me think that the second option
might be truly fruitful only if pursued as a complement of the third. It is, then, in the frame of this third option that I reason here.

When this option is joined with what I have  said at the end of \S~\ref{sect.PS}, above, it brings me to consider that a
crucial task for the philosophy of mathematical practice is that of suggesting good
conceptual categories to be used to account for the structural form of mathematics (intended as an existent reality), as it appears to be.

I say `account for the structural form' for the sake of brevity, of course, in order to make clear that the relevant conceptual categories are not those required to
describe the specific contents of mathematical theories, but rather some general features of these theories and contents, which are often characterized, by using
a technical philosophical language, as logic, epistemic, ontological, etc. The specific  contents of mathematical theories
are also to be described, of course, using appropriate conceptual categories, but shaping these categories, and providing such a description is a task for 
mathematicians rather than philosophers  (though it would be wise for philosophers to try to understand this description, without falling into the opposite temptation
to that mentioned in \S~\ref{Sect.Mis}: the temptation of replacing mathematicians in this task, to provide a supposedly deeper description
of these contents, which inevitably results much less profond than confused 
and defective; a temptation, it should be said, that some alleged philosophers of mathematical
practice are not unfamiliar to).

In this framework, it seems to me that nothing might forbid philosophy of mathematical practice from embracing a platonist perspective,
under the form of the thesis that the peculiar nature of mathematical activity is to (or, at least, might suitably) be 
accounted for as an intellectual activity dealing with (the) abstract objects (denoted by the usual terms figuring in the mathematical language).

This is still a very broad view, however, that can be declined in very different ways. One might wonder, for example, what one should mean in this context
by `intellectual activity', `dealing with', `abstract objects', and, even, by `objects' \textit{tout court}. I cannot enter all these matters, here, of course. I will, then, take
for granted that my audience understands the term `intellectual activity' more or less as I do. I only add that, according to the version of platonism that I'd
like to illustrate (and defend) here, `dealing with $x$' means something like `supplying (shaping, giving birth to) $x$, 
then studying it in order to describe it'. I also suppose, that
this is enough to restrict my ongoing task to (quite programmatically) answer to the following single question: what should one mean, in the present context, by
`objects', in general, and `abstract objects', in particular. The attention to the context is essential, since I'm far from pretending to explain what can make
an object abstract, in general. I shall limit myself to the case of abstract objects that one might take as mathematical ones. 

My purpose, in what follows, will, then be nothing but that of trying to
make clear as the notion of an abstract (mathematical) object should be understood, in order to provide a good conceptual category to be used while describing
mathematics as an intellectual activity dealing with abstract objects.

\section{Objects by Reification}

Objects are usually taken to be self-existing individuals, or, better, individual self-existing
contents. An object is, then, taken to be abstract if this content is so, that is, as it is usually said, it is not spatiotemporally identifiable or determinable. An 
abstract object is, then, taken to be mathematical if (pure) mathematics involves putative reference to it (without receiving it from outside, for instance from logic).
If we understand an object this way, it becomes unavoidable to wonder how an individual content can exist (or have existed, or will possibly exist) without 
being spatiotemporally identifiable or determinable. It is just by taking this possibility for granted, without any further explanation, or, worse, in trying to
provide quite improbable, and, often, even (quasi-)mystical answers to this question that platonism becomes generally implausible.

A way to avoid the problem is by ceasing to take existence, at least if intended as a primitive notion, as a necessary condition for objecthood. This is in no way
an original perspective, and it is rather quite common in philosophy, also outside philosophy of mathematics. It goes, however, quite often together
with some troubles, which is neither necessary nor possible to consider here in any detail. It will be enough to observe that a common move made to escape (some
of) these troubles is taking abstract objects to be individual contents coming from reification of (first-level) concepts. This could go together with making the notion of an 
abstract object essentially independent of the more general notion of an object, namely with renouncing at defining abstract objects as  a particular sort of objects.
What is important here is, however, not how this last task could be achieved, but rather how one can make the notion of reification of a (first-level) concept both precise 
enough, and such to avoid that the claim that mathematics deals with abstract objects becomes empty.

An example will make the point clearer. Let us take for granted that ordinary arithmetic deals with natural numbers. A very  low-cost way to also grant that
this makes it deal with (abstract) objects is by providing a clear enough 
definition of the concept of a natural number, by means, for instance, of some appropriate version of Peano axioms, then taking natural numbers to 
be abstract objects just insofar as they are reifications of the concepts of the distinct items falling under this concept, reiteratively defined in turn, on the base
of the definition of this last concept. The notion of reification of a concept is not very well defined here, but it does not 
seem \textit{prima facie} implausible to consider, for instance, that 
taking the number one to be the successor of zero, by also admitting that the concepts of zero and successor are implicitly defined 
through the relevant version of the Peano axioms, be
a way to take this number as the reification of a concept (that of successor of zero), and that the same also happens, 
in principle, for all other positive natural numbers.

Though certainly not perfectly crystal-clear, this reasoning does not seem to me preposterous. But it would not
help up to make the claim it supports---that arithmetic deals with abstract objects, understood as reification of concepts---anything more than an empty 
(though comfortable) claim, a contrived  \textit{fa\c{c}on de parler} to merely say that mathematics uses singular terms with a referential intent.
If we want to make platonism an interesting philosophical thesis, we need to do much better than this.

Today's philosophy of mathematics presents at least two ways of doing it, each of which has unquestionable merits. 

The first is provided by \textit{Ante Rem} Structuralism (ARS), mainly advocated by Steward Shapiro (\citeyear{Shapiro1997};
a similar view, that differs, however, from Shapiro's in many important details is advanced in \citealp{Resnik1997};
other forms of structuralism are also available, which is not necessary to discuss here, since they do not openly adhere
to the idea that mathematical objects result from a form of reification of concepts). 
To say it in general, the basic idea behind ARS 
is characterizing the appropriate
form of reification by identifying a particular way, typical of mathematics, for defining concepts and organize them with the
aim of fixing individual contents. This modality pertains to the identification of a categorical structure: the relevant abstract
objects are identified with the poles (or places) in such a structure.  They are, then, taken to exist insofar as the categorical
structure they are poles of exists, and this exists, in turn, if and only if it is coherently defined, and obeys a
general theory of structures, which (\textit{mutatis mutandis}) mimics set theory.  

Though much more might be said on ARS, this short account should be enough for making
clear the way it answers to Platos's problem. What is important, indeed, is that it leaves to mathematics
itself the task of defining 
categorical structures, and then insuring 
then existence of mathematical objects, though taking on it the task of indicating the appropriate way
that has to be done (namely by defining categorical structures). Leaving apart the
usual, more technical objection about indiscernibles (originally advanced by \citealp{Keranen2001}, to which Shapiro
has replied, for instance, in \citealp{Shapiro2008}), its main shortcoming
comes, to my mind, from its being unapt to account for the historical evolution of mathematics, and, more in general,  for large portions
of past and present mathematics which 
are certainly not structuralist in nature, as well as for the current pervasive practice of transcendental proofs (proofs of theorems of a certain theory, 
involving tools from other theories).

The second way requires a more detailed explanation, since it appeal to less commonplace technicalities. It also adopts 
a much more descriptive (rather than normative) stance. It is provided by E. Zalta's Object Theory (\citealp{Zalta1983}, \citeyear{Zalta2021}). 
The basic idea of this theory (or OT from now on)  is  providing a
logical rendering of the very idea of reification of concepts in the shape of a particular form of predication, called `encoding', essentially distinct from that 
customarily in use in predicate logic. Formally, this depends on admitting two distinct forms for an atomic predicate formula, namely `$Fa$', and 
`$aF$', where `$a$' is a term and `$F$' a monadic predicate, expressing the case where $a$ respectively exemplifies (or has) $F$,
and encodes it. Both forms of predication are embodied in a system of modal higher-order logic admitting types, but also working without them for simplicity.
In its non-typed version, this system includes the following axiom
$$
\Diamond \left[aF\right] \Rightarrow \square \left[ aF\right]
$$ 
(where `$a$' and `$F$' are schematic letters, and `$\Diamond$' and `$\square $' are the usual modal operators), asserting that $a$ can encode $F$ only if it encodes it necessarily. As the reciprocal implication also holds, of course, because of the usual properties of the modal operators, this renders that idea that $a$ is the reification of $F$ if and only if it is necessarily so.

This idea goes together with two further crucial ingredients of the system: the clause that an object encodes a (monadic) property only if it is abstract,
and a form of comprehension asserting that, for any condition on monadic properties expressible in the formal language at issue, there is an
abstract object that just encodes the properties that meet this condition
(which allows an abstract objet to be conceived, in fact, as the reification of a bunch of properties).

The former ingredient is implemented by the introduction of the non-logical monadic predicate constant `$E!$' denoting the property of being a concrete object,
through which the two other monadic primitive predicate constants `$A!$' and `$O!$', respectively denoting the properties of being an abstract and an ordinary one, are defined by
$$
A! =_{df} \left[\lambda x \lnot \Diamond E!x\right] \qquad \text{and} \qquad O! =_{df} \left[\lambda x \Diamond E!x\right].
$$
It follows that any object considered in the system is either abstract or ordinary, and it
is abstract if and only if  it cannot be concrete, that is, it is concrete in no possible world, while it is ordinary if and only if it can be concrete, namely
it is concrete in some possible world. To these definitions, the following axiom-schema is added:
$$
\forall \ldots \left[O!x \Rightarrow  \lnot \exists P \left[xP\right]\right],
$$
(where `$\forall \ldots \left[\varphi \right]$' stands for the universal closure of $\varphi$),
asserting that an objet is ordinary (and, then, not abstract) only if it fails to encode (is not the reification of)
any (monadic) property. By contraposition, it immediately follows that an object encodes (is the reification of)
of a (monadic) property only if it abstract, as announced above.

The latter ingredient completes the picture. It consists of the following comprehension axiom-schema:
$$
\exists x \left[A!x  \wedge \forall P \left[xP \Leftrightarrow \varphi \right] \right]
$$
where `$\varphi$' is any formula in the formal language at issue in which `$x$' does not occur free.
Let, then, $\boldsymbol P$ be whatever (monadic) property already defined in the system. From this
axiom schema, together with the first axiom mentioned above and a principle of identity for abstract objects
making $x$ the same abstract objects as $y$ if and only if $x$ necessarily encodes the same properties as $y$,
we immediately get that
$$
\exists ! x \left[A!x  \wedge \forall P \left[xP \Leftrightarrow  \forall z \left[Pz \Leftrightarrow \boldsymbol Pz \right] \right] \right],
$$
which ensures that there exist a single abstract object that encodes (is the reification of) any (monadic) property coextensive
with  $\boldsymbol P$. Extensionally speaking, we can then say that for whatever (monadic) property $\boldsymbol P$, there is a single
abstract objet that encodes (reifies) it and only it.  But this is not all, since the formula `$\varphi$' can be more complex than that.
It can express any condition on (monadic) properties, which makes the schema ensure that for any such condition there is an
abstract object that just encodes the (monadic) properties that meet it, as also announced.

With this formal apparatus at hand, and with the help of some further enrichments of the language of his 
theory---including the introduction of the definite descriptions operator `$\iota$' (such that `$\iota x \left[\varphi\right]$' is a term
denoting the single object, if any, that meets the condition expressed by `$\varphi$')---,
B. Linsky and Zalta himself have been able to provide a rigorous explicit characterization
of a intra-theoretic mathematical object (\citealp{LinskiZalta1995}, \citeyear{LinskiZalta2006}, \citeyear{LinskiZalta2019};
the reason for I speak of an  intra-theoretic mathematical object, rather than of a  mathematical object 
\textit{tout court} will become clear below). 
This is an object $\kappa_{\textsc{t}}$ such that there is a mathematical theory $\textsc{t}$ (identified, in turn, with an
abstract object appropriately fixed within OT, so as to render an actual mathematical theory), such that this object encodes all and only the
properties it has in ${\textsc{t}}$, that is:
$$
\kappa_{\textsc{t}} = \iota x \left[A!x \wedge \forall F \left[xF \Leftrightarrow \textsc{t}\left[\lambda z F \kappa_{\textsc{t}}\right] \right] \right],
$$
where $\left[\lambda z F \kappa_{\textsc{t}}\right]$ is the $0$-place (or propositional) 
property that any object has just in case
the world it belongs to is such that $F \kappa_{\textsc{t}}$ (in Zalta's terminology, this is the propositional property ``constructed out'' of $F \kappa_{\textsc{t}}$).

This is properly not a(n explicit) definition, of course. Its very logical form prevents it from being so, since the term `$\kappa_{\textsc{t}}$' figures in it 
in the left-  and the right-hand sides. It is rather a ``theoretical
description'', offered in OT's (rigorous) language. Whether there is 
a theory like $\textsc{t}$, if it is actually a mathematical one, how it its done, 
which individual contents (denoted by singular terms) it deals with, which of them are taken to have in it 
this or that  property, and what makes they have them are questions that OT neither can decide nor pretends to decide. 
It simply provides a formal 
setting making possible to rigorously claim that these contents are intra-theoretical 
mathematical objects, and that nothing else is
such an object. This entails that intra-theoretical mathematical objects are all abstract, 
and are just the individual contents (possibly denoted by appropriate singular terms) that a 
mathematical theory deals with. More generally, each of them is the reifications of
the concept of having all the properties that this theory assigns to it. 

Just as Shapiro's structuralist
definition, Zalta's characterization or ``theoretical description'' of intra-theoretic mathematical objects is intrinsically unsuited 
for accounting transcendent proofs and extra-theoretical contents, that is, for mathematical objects possibly studied by 
different theories and then independent of them. The problem is serious, since 
denying this possibility makes it impossible to maintain, for instance, if not quite informally, that different versions of arithmetic 
all deal with natural numbers, as a platonist account should be able to do, it seems to me.
OT has however the ressources for also characterizing extra-theoretical mathematical objects.

On the one side, it accounts for the reference of singular mathematical terms in statements that make no
(explicit) mention of a particular mathematical theory, as a reference to abstract objects that encode a relevant
concept. Supposing we are able to deal with the concept of a triangle, however it might be
fixed, it would be enough, for instance, to identify this concept with a particular concept, call it $\boldsymbol T$,
dealt with within this theory, to easily render the statement of a natural or informal language $\ulcorner a$ is a triangle$\urcorner$ 
by the statement  of such a theory $\ulcorner a\boldsymbol T\urcorner $. 
This is the same as interpreting the `is' in the former  statement as an informal occurrence of the encoding operator.
If it happens that the concept  $\boldsymbol T$ is defined within a certain particular mathematical theory, be it formal or not
(and this theory is, in turn, identified with an abstract object appropriately 
fixed within OT), 
we can prove (within this theory) that $a$ encodes $\boldsymbol T$ if and only if $a$ exemplifies $\boldsymbol T$. The statement 
$\ulcorner a\boldsymbol T\urcorner $ is then a statement about an intra-theoretic mathematical object,
and this makes it  suitable for rendering $\ulcorner a$ is a triangle$\urcorner $ only if the same also happens for this
latter statement. Still nothing forbids considering $\boldsymbol T$ (informally) fixed otherwise (by a way that
OT is not required to specify), and still rendering $\ulcorner a$ is a triangle$\urcorner $ by  $\ulcorner a\boldsymbol T\urcorner $.

On the other side, OT also allows to define within it a quite weak, pseudo-freegan set theory, and use it to
define natural number \textit{\`{a} la} Frege, as extensions of concepts \citep{Zalta1999}.  To this purpose,
it is enough to take the extension $\varepsilon \boldsymbol F$ of a (sortal) concept $\boldsymbol F$ as the object
$$
\iota x \left[A!x  \wedge \forall P \left[xP \Leftrightarrow  \forall z \left[Pz \Leftrightarrow \boldsymbol Fz \right] \right] \right]
$$
that just encodes $\boldsymbol F$, to define the appurtenance relation $\in$ (or its Fregean \textit{alias} $\frown$)
by
$$
x \in y \Leftrightarrow \exists P \left[x=\varepsilon P \wedge Py\right]
$$
and to adopt an appropriate modal principle to get second-order Peano axioms as Frege did (\textit{mutatis mutandis}) 
in \textit{Grungesetze}. According to Zalta, extra-theoretical natural numbers can, then, be taken as the objects meeting
these axioms (on which no order relation, no addition, and no multiplication is defined yet), and the same can be done for
other mathematical objects definable within the same weak set theory.

All this should be enough to make clear, in the same time, the great expressive (or descriptive and elucidatory) strength of OT, 
and its intrinsic and purposeful presumptive weakness. This weakness consists in its deliberately avoiding,
any endorsement in favor of a
specific view about the nature of mathematics. The theory succeeds in providing a clear logical account of the reification process
and the corresponding  constitution of abstract objects. 
However, it purposely refrains from accounting
for what distinguishes the mathematical form of reification, and makes the intellectual exercise that
mathematics pertains to result in such a powerful and stable corpus. This makes OT particularly suited
as a theoretical context for advancing a pluralistic view about mathematics \citep{ZaltaFC}. But it makes it also
unsuited to provide, alone, an answer to Plato's problem. This does not mean that it cannot contribute to 
provide this answer. It simply means that some other independent insight is required to achieve this task.

\section{Object by \textit{De Re} Epistemic Access}

This is why I want to suggest  another perspective here, essentially independent  of ARS and OT.
In particular, I suggest looking at an object, both concrete or abstract, as
an individual content we could have (had) a \textit{de re} epistemic access to. 

If I say $\ulcorner a$ has the property $\boldsymbol P\urcorner$ I can mean two quite different things, whose difference depends on the modality of my epistemic access to $a$.
I can mean, of $a$, identified as such in advance, that it has the property $ \boldsymbol P$, or I can mean that $a$, 
which has not been necessarily identified as such in advance, 
has the property $\boldsymbol P$. 

I was born in Varese, a nice little town in the north of Italy, and lived permanently there the first twenty-seven years of my life 
and, even now, though I live
abroad, I try to come back to Varese anytime I can, spending time in the house where I grew up and lived for so many years. 
From the balcony of this house, one can see the
Varese Lake, on the shores of which I went to walk and bike myriads of times. I also went boating many times on this lake, and, when I was a kid, pollution
did not advise against swimming in it, and so I swam in it many times. 
Now, it just happened that I read on Wikipedia that the Varese Lake has a surface area of 14.5
km$^2$. Despite my familiarity with Varese Lake, this is a new piece of information for me, 
and it is, literally speaking, a piece of information about Varese Lake. It is about
it, that I got this information, which I did not have before (simply since I never thought before to look for it). 
Not having this information did in no way made me unable to perfectly
know what the Varese Lake is, to understand what one was speaking about when speaking of it, 
to locate it on the surface of the Earth, and the image of it in my
memory: in few words, to have a perfect epistemic access to it. Because 
of this access that I now know what the new information that I  got is about. So, if I 
now say,  $\ulcorner$The Varese Lake has a surface area of  14.5 km$^2\urcorner$, I definitely mean, of the Varese Lake, 
that it has a surface area of 14.5 km$^2$.

After having got this new information about the Varese Lake, I opened another Wikipedia page, and did read there that  the 
``Eridania Lake is a theorized ancient lake 
[on Mars] with a surface area of roughly 1.1 million square kilometers [\ldots][,~a] maximum depth of 2,400 meters [and a] volume[\ldots][of] 562,000 km$^3$''. 
I had never heard anything concerning Eridania Lake, and never heard of the theory that there could have been lakes on Mars in the far past.  
I'm, then, very much surprised that it could have existed such a huge mass of water on Mars. This new piece of information is not about something I knew 
before having it and independently of it, that I was able to locate in the universe, and I had an image of in my memory. In few words, it is not about something that
I had an epistemic access to, before getting it. This piece of information is instead, for me, constitutive, of the Eridania Lake: it is what makes me have
an epistemic acces to it.
So, if I now say $\ulcorner$The Eridania Lake had a surface area 
of  about 1,100,000 km$^2\urcorner$, I cannot but mean, literally speaking, 
that the Eridania Lake had such a surface area, and this is what  makes me now able to think at it somehow.

This illustrates the distinction I want to make between two radically different modalities of epistemic access to a content that works as a reference 
of a proper name. In the former example, I take my epistemic access to be \textit{de re}. In the latter, I take it to be purely  \textit{de dicto}. 

This having been said, let us come back to the latter example (the former seems to me clear enough for not requiring other considerations).  Though no 
human being (and also no nonhuman animal presumably) has certainly never had any direct empirical experience
of the Eridania Lake,  it seems to me clear that it counts as an object, and even more so as a concrete one. What does make it so clear (at least to me)? 
I advance that the reason is this: the fact that nobody has ever had any direct empirical experience of the Eridania Lake does in no way make it
logically impossible for someone to have had such a direct empirical experience. Moreover, even if 
(for some reasons, that I am not frankly able to imagine) we were urged to consider
such an experience as  logically impossible, it would still remain that this would hardly depend on the intrinsic features of the Eridania Lake, 
but rather of the limitations imposed to our (or other animal's) capacity of empirical experience. As such, the Eridania Lake would have 
allowed anyone to have a direct experience of it, and this experience
would have, in turn, allowed a human  being to have a \textit{de re} epistemic access to it.
The (logical) possibility of such an access is, in my view, what makes the Eridania Lake count for us as an object. The fact that this is the 
possibility of an empirical experience is, instead, what makes this object count for us as a concrete one.

There is an important aspect of this argument that deserves clarification.
What I just described are in no way conditions that
make the Varese and the Eridania Lakes concrete objects \textit{tout court}, 
but rather conditions that make them count for us as concrete objects. I am, in fact, not interested in
the question of establishing what makes a concrete or an abstract object be so. I even consider this question perfectly meaningless. 
In agreement with the approach to philosophy 
that I described above, I'm rather interested in the question of establishing what makes us consider an object to be so. It is just to this question that
I suggest answering by saying that this is nothing but
the possibility of having (had) a \textit{de re} epistemic access to it.

This having been made clear, let us come back to the main point. The two examples I advanced are both such to make the possibility of having (had) 
a \textit{de re} epistemic access to a certain individual content depending on the possibility of a direct empirical experience of this content. This
is just what makes make these examples concern concrete objects.  The next question is, then, this: are there other sorts of epistemic relations with some 
individual contents that can make us able to have a \textit{de re} epistemic access to them. The possibility of (sensibly) taking 
some individual content as an abstract
object---rather than, merely, as the putative referential target of a proper name used in a statement that is admitted to hold---depends, 
in my view, on the answer to this
question. Moreover, if the answer is positive, it is just the nature of 
these alternative relations that decides  the nature we can ascribe to the object at
issue, and whether it should count for us as an abstract object or a concrete one. 
This makes me replace the question of whether there are mathematical (abstract) objects with the question of whether there are forms of epistemic relations
with individual (abstract) contents that are proper to mathematics (and differ, then, from direct empirical experience) that make 
it possible for us to have a \textit{de re} epistemic 
access to these contents.

To better understand the question, let me try to better clarify it. 

What makes my epistemic access to the Varese Lake \textit{de re} is my capacity of identifying this content as such 
while ascribing a certain property to it, but fully independently of my doing it. This is, in other terms, my capacity of fixing this 
content in my epistemic horizon independently of recognizing it as the content having this property.

Now, it seems plain to me that this also happens very often with individual contents that we deal with in
mathematics. 

We can certainly identify the number 31 independently of ascribing to it the property of being the eleventh prime number
(in the usual order on natural numbers).
Hence, properly speaking, we do not say that  31 is the eleventh prime number; we rather we say, of 31, that it is the eleventh prime number. 

In the same way, 
we do not properly say that the ratio of two successive terms of the Fibonacci sequence tends to the golden ratio, or even that this ratio is (equal to) 
$\displaystyle \frac{1+\sqrt{5}}{2}$; we say of the ratio of two successive terms of the Fibonacci sequence, that it tends to the golden ratio, and,
of this ratio, that it is equal to $\displaystyle \frac{1+\sqrt{5}}{2}$. 

Despite naming it after Euclid, we do not properly say that the 
Euclidean closure of the rational numbers is the field of the numbers constructible by ruler 
and compass;  we say, of the Euclidean closure of the rational numbers, that it is the field of the  numbers constructible by ruler 
and compass. 

Again, we do not properly say that $4$ is the minimal number of colors required to color a map by avoiding 
that two adjacent regions be colored with 
the same color; we say, of $4$, that it is the minimal number of colors required for so coloring a map. 

We no more properly say that $18$ is 
the Ramsey  number $R\left(4,4\right)$ (i.e. the smallest number of vertices of a graph containing either a 4-cycle or 
an independent set with $4$ vertices), or
that $43 \leq R\left(5,5\right) \leq 48$; we say, of $18$, that $R\left(4,4\right)$ is equal to it, and, of $43$ and $48$, that they are the limits of the
closed interval including $R\left(5,5\right)$.

Finally, we do not properly say that $\displaystyle{\frac{\pi^2}{6}}$ is
the limit of the series $\displaystyle{\sum_{i=0}^{i=\infty}\frac{1}{i^2}}$; we say of  $\displaystyle{\frac{\pi^2}{6}}$, that it is the limit of this series. 

Examples might be multiplied. Some of the previous ones might, by the way, also be inverted. We might perfectly say, for instance, of 
$\displaystyle \frac{1+\sqrt{5}}{2}$, that it is (equal to) the golden ratio, or of the limit of the series $\displaystyle{\sum_{i=0}^{i=\infty}\frac{1}{i^2}}$, that
it is (equal to) $\displaystyle{\frac{\pi^2}{6}}$. The point is not which examples we consider, or what we consider to be, in each of them, the given 
content, and the property that we attribute to it. What is relevant is that these examples illustrate an important fact about mathematics: that in it,
we can exhaustively fix individual contents in our epistemic horizon, independently of ascribing to them relevant
properties that are rather ascribed to them \textit{post festum}. In other terms: we can identify such contents without 
specifying all the properties they have in the very setting in 
which they are identified or in other settings. This fact is just what allows us to export individual contents from setting to setting, or, better, to deal with the same 
individual contents in different settings.  Again, we can fix such contents in a certain setting, 
then extend the setting by defining appropriate relations or functions 
on these same contents, which allows us to ascribe new properties to them, without making them different contents. We do both things 
pervasively. And this is just what makes, to my mind, a sort of spontaneous (though often irreflexive) 
platonism so popular among working mathematicians. 

This is, indeed, my own version of platonism for a philosophy of mathematical practice: mathematics can be (is) suitably accounted for as an 
intellectual activity  dealing with \textit{sui generis} individual abstract contents, so fixed as to make us able to have 
a \textit{de re} epistemic access to them. As it should be clear, this claim has strictly nothing to do with the existence of these contents 
(and the truth of some statements on them). Not only does it not assert the existence of mathematical objects, it also
does not deny it. It rather assumes that the notion of existence of abstract objects is, when conceived as a primitive notion, simply
nonsensical, by leaving open the possibility of defining this notion with base on other ones, by so ending up recognizing that, yes,
mathematical objects exist (and mathematical theorems are true of them), but in a new, derived sense.
Among other things---particularly the admission that mathematical statements are to understood \textit{prima facie}---,
this assumption of nonsensicality essentially distinguishes the claim from any version of nominalism. It is, rather a claim for genuine platonism, though
this is,  as it were, a platonism without existence (and truth), at least if intended as primitive notion(s).

I say `sui generis', since the way the relevant individual contents are fixed is peculiar to mathematics. And it is just this peculiarity that makes this thesis not only 
compatible with the program of the philosophy of mathematical practice, but also dependent on this very program for its clarification and advocacy. Indeed, it is only
by parsing mathematical practice and its history that we can identify the different ways these contents 
are fixed, so as to make such an epistemic access possible, and we can
try to recognize what makes them mathematical in nature. 

Such a project is structurally similar, 
in nature, to that of any empirical science: it is, so to say, the project of
an empirical science of mathematical reality, this reality being not, however, that of self-standing objects, as a traditional platonist would argue for, but
rather that of mathematical activity and production.

\section{History, Cognition and Logic}

Of course, realizing this project is not something I can do here. Both the limits of space and the very purpose of my essay force
it to be essentially programmatic, in nature. What I can and want to do is delineating a philosophical perspective,
a research program for a philosophy of mathematical practice able to answer the Plato's problem. This only
requires, for the time being,
to briefly suggest three general sources (certainly not alternative, but rather complementary) 
for our ability to fix  abstract individual contents so as to make possible for us to have \textit{de re} epistemic access to them
and, then, to take them to be mathematical objects): history, cognition, and logic (or, more generally, formalization).

\subsection{By History}  

Mathematics largely evolves by (conceptual) recasting (I borrow this term from several talks of Ken Manders). Theories are firstly (mainly informally)
advanced, then variously restated, reorganized,  formalized. This important part of mathematical activity explicitly and essentially depends on the admission that
these different forms of recasting conserve contents, that some contents (not only individual, of course, but crucially also individual ones) remain the same
under the transformation of the setting in which they are defined and studied. This
suggests to inverse the model-theoretic perspective. 

According to this last perspective, mathematical objects are conceived of as the elements of the 
first-order domain(s) involved in the model(s) of a (formal) mathematical theory, namely with the individuals forming the range(s) of the (different sorts of)
first-order variables. 
What I mean by inverting this perspective is conceiving most mathematical theories as codifications---or, possibly, in some cases, formalizations---of a portion 
of available knowledge inherited from the history of mathematics (in the proper sense of `history', i.e. as \textit{res gest\ae}, rather than \textit{historia rerum gestarum}), in which contents (not only individual, but crucially also individual ones) are 
crystallized, so as to make their transmission possible. 
What many mathematical theories would, then, do is taking some contents as given, 
and suggest a new conceptual configuration
for them, possibly connecting them with other contents, historically crystallized in turn, but dissociated from the former ones, in the original piece of knowledge. 
This would make a mathematical theory be not about its models but about the contents it captures by recasting. 

The availability of a single model, under
isomorphism (and within a given semantic setting) would merely show that the theory captures these contents making them part
of an univocal structure (within this setting). The presence of non-isomorphic models would rather show that these contents 
are not so captured, leaving, then, space
for alternative interpretations. This lack of univocality (or relative categoricity) can be certainly beneficial, since it can open new horizons of knowledge (the case of 
non-standard analysis is symptomatic), but it says nothing about the original contents. It rather says something about new contents that, in this way, 
we become able to distinguish from the former, and study together with, or separately from them.

Here, the basic idea is that mathematical activity is intrinsically historical (as argued for in the preamble), 
though this characteristic and peculiar nature of it remains hidden to superficial regards. The reason of the concealment is 
that this historical nature manifests itself by a sort of (temporary) negation of itself, just due to recasting.

\subsection{By Cognition}

\noindent Mathematics requires (both common and specific) cognitive abilities. These are not the usual abilities to realize practical performances, 
often considered by alleged studies in mathematical (typically numerical) cognition, in order to display widespread forms of cognitive access 
to mathematical objects
(typically small integer positive numbers), some of which would also be  proper 
to human infants and some species of primates. These merely appear to me as mundane abilities to 
play some linguistic or prelinguistic games, which either involve
mathematical words, or admit a retrospective description appealing to these words or the corresponding concepts. 
Subitizing is only one among many other 
possible examples of these practical performances, which have, to my mind, no attested relation with mathematical theories and the 
contents they actually deal with.

What I'm instead referring to are cognitive abilities that are required to crystallize intellectual contents involved in abstract reasoning. 
These are typically (but not limited to) abilities of recognizing an invariant structure within different (sorts) of exemplification of it, for example the 
$\omega$-sequence structure within the succession of Arabic numerals, or that of the reciprocal of the natural numbers, or the structure of a regular polygon
within equilateral triangles,  squares, pentagons, etc. 

More generally, these are abilities of detecting abstract invariants, either individual or not (for instance 
operational ones) and of both considering them as such, and working on particular instances of them, in order to identify general structural properties. 
In Kant's well-known terminology, this is the capacity of dealing
with the ``universal \textit{in concreto}''. However flawed his conclusions may have been, both by a dramatically mistaken experimental practice and by a
misleading confusion between ontogenetic and historical evolution, in his study of 
the cognitive development of children, J. Piaget famously tried to locate and isolate abilities like those

Here, the basic idea is that performing mathematics activity depends of the (acquisition and) possession of cognitive abilities of different levels and sorts.
These are the abilities that make us
able to identify stable intellectual contents (either individual or not), which are both phylogenetically forged, and ontogenetically 
refined, also through cognitively complex and specific exercises (even if certainly not only through them).

\subsection{By Logic (or Formalization)}

\noindent Mathematical (conceptual) recasting mostly depends on a codification of 
language and intellectual acts, such as (stable or provisional)
suppositions, inferences, definitions, constructions, etc. This codification is essentially logical, insofar as it depends on restrictions of the forms of expression and
of the possibility of actions, in agreement with a number of general clauses to be applied in any specific occasion. These clauses might be seen as norms 
granting permissions, to be appealed to yield an argumentative setting that excludes all that is not allowed by them. 
A given setting is, of course, often extended,
either occasionally, by unenvisaged admissions or acts, or systematically, by the addition of clauses granting further permissions. Moreover, different settings are
frequently merged with each
another. If  Cantor was certainly right in claiming that ``the essence of mathematics just lies in its freedom'' \cite[p.~564; my translation]{Cantor1883}, this is essentially the freedom of fixing appropriate 
settings or, at
most, of locally contradicting them by merely foretelling new permissions, rather than that of leaving anything go. 

So broadly described, codification can either result in
formal or informal theories. The border between the two cases is far from rigid and largely depends on historical or cultural conditions,
branches of study, and methodological attitudes and sensibilities. What can be taken as formal within some mathematical communities can not be taken this way within 
other communities. Even the criterion of (the possibility of a) computer implementation is much vaguer than it might appear at first glance, since it depends on the machine
languages used and the interfaces with higher-level languages. Still, be it as it may, the crucial point is that any form of codification, as broad as it might be,
goes together with  identifying a delimited space of possibilities, involving some poles that can count as fixed individual contents. 

Conceiving a particular
codification as an act of recasting makes these individual contents count as particular occurrences of a more general content that they render within a certain setting.
This is just what allows the last content to be taken, as such, as an available target for attributions of properties or the definition of relations or functions, which are in no way
constitutive of it, either in general or within the particular setting at issue. 

Fixing natural numbers as poles of the (relatively) categorical structure defined by Peano
second-order axioms allows, for instance, both taking these poles as targets of a fist-order (uncategorical) axiomatization---which count, then, as an act of recasting
within a different setting---, and defining on them, within their own constitutive setting, a linear order, and an additive and a multiplicative law of composition,
without changing their nature, but rather attributing to them, as such, a new bunch of properties.
Moreover, making codification meet some appropriate conditions also allows its outcome to become, as such, the subject matter of a new mathematical
inquiry, depending, in turn, on a meta-codification. Such is the very well-known phenomenon of meta-mathematics, which is also a distinctive aspect of our capacity
of having, in mathematics, a \textit{de re} epistemic access to individual contents.

Here, the basic idea  is that mathematics uses codified languages within delimited argumentative spaces. These spaces are not only suitable for fixing some
appropriate structures that are intended to render previously available pieces of knowledge, by  bringing them to a new life. They also have specific structural properties that can be studied  as such by new mathematical theories.

\section{A Few Words, in Guise of Conclusion}
\noindent  When they are so conceived, as sources of our capacity to fix individual contents, so as to make possible having a \textit{de re} epistemic access to them,
history, cognition and logic become crucial ingredients of a form of platonism. Its being independent of any 
talk of existence of abstract objects (and of truth of statements on them), if not appropriately re-defined as derivative notions, makes the theses advanced by it perfectly factual claims. 

Whether mathematics complies or not with such a form
of platonism, becomes, then, a pure matter of fact: a matter of fact that the philosophy of mathematical practice is not required to abstractly argue for, but rather to
verify by an inquiry of this practice and its history. This makes this philosophy work, in turn, as a sort of high-level empirical science.

\bibliographystyle{plainnat}
\bibliography{PMP.bib}

\end{document}